\documentclass[10pt]{article}
\usepackage{latexsym, amsfonts, amsmath, amssymb, amscd, epsfig}

\numberwithin{equation}{section}
\newtheorem{thm}{Theorem}[section]

\newtheorem{lem}[thm]{Lemma}

\newtheorem{re}{Remark}[section]
\newenvironment{pf}{{\noindent \it \bf Proof:}}{{\hfill$\Box$}\\}

\begin{document}

\title{A regularity criterion of the 3D MHD equations involving one velocity
and one current density component in Lorentz space}
\author{Ravi P. Agarwal$^{1}$, Sadek Gala$^{2}$ and Maria Alessandra Ragusa$%
^{2,3}$ \\
$^{1}${\small Department of Mathematics, Texas A\&M University-Kingsville,
Kingsville, USA, }\\
{\small Ravi.Agarwal@tamuk.edu}\\
$^{2}${\small Dipartimento di Matematica e Informatica, Universit\`{a} di
Catania,}\\
{\small Viale Andrea Doria, 6 95125 Catania - Italy,}\\
{\small sgala793@gmail.com}\\
$^{3}${\small RUDN University, 6 Miklukho - Maklay St, Moscow, 117198,
Russia,}\\
{\small maragusa@dmi.unict.it}}
\date{}
\maketitle

\begin{abstract}
In this paper, we study the regularity criterion of weak solutions to the
three-dimensional (3D) MHD equations. It is proved that the solution $(u,b)$
becomes regular provided that one velocity and one current density component
of the solution satisfy%
\begin{equation}
u_{3}\in L^{\frac{30\alpha }{7\alpha -45}}\left( 0,T;L^{\alpha ,\infty
}\left( \mathbb{R}^{3}\right) \right) \text{ \ \ \ with \ \ }\frac{45}{7}%
\leq \alpha \leq \infty ,  \label{eq01}
\end{equation}%
and
\begin{equation}
j_{3}\in L^{\frac{2\beta }{2\beta -3}}\left( 0,T;L^{\beta ,\infty }\left(
\mathbb{R}^{3}\right) \right) \text{ \ \ \ with \ \ }\frac{3}{2}\leq \beta
\leq \infty ,  \label{eq02}
\end{equation}%
which generalize some known results.
\end{abstract}

\noindent {\textit{Mathematics Subject Classification(2000):}\thinspace
\thinspace 35Q35,\thinspace \thinspace 35B65,\thinspace \thinspace 76D05.}
\newline
Key words: MHD equations; regularity of weak solutions; Lorentz spaces.
\newline

\newpage

\section{Introduction}

\bigskip This paper deals with the well-known problem of the regularity of
the solutions for the 3D magneto-hydrodynamical (MHD) system
\begin{equation}
\left\{
\begin{array}{c}
\partial _{t}u+u\cdot \nabla u-\left( b\cdot \nabla \right) b-\Delta
u+\nabla \pi =0, \\
\partial _{t}b+u\cdot \nabla b-b\cdot \nabla u-\Delta b=0, \\
\nabla \cdot u=\nabla \cdot b=0, \\
u(x,0)=u_{0}(x),\text{ \ }b(x,0)=b_{0}(x),%
\end{array}%
\right.  \label{eq1.1}
\end{equation}%
where $u=(u_{1},u_{2},u_{3})$ is the velocity field, $b=(b_{1},b_{2},b_{3})$
is the magnetic field, and $\pi $ is the scalar pressure, while $u_{0}$ and $%
b_{0}$ are the corresponding initial data satisfying $\nabla \cdot
u_{0}=\nabla \cdot b_{0}=0$ in the sense of distribution.

Since Duvaut-Lions \cite{DL} and Sermange-Temam \cite{ST} constructed the
so-called well-known weak solution $(u,b)(x,t)$ of the incompressible MHD
equation for arbitrary $(u_{0},b_{0})\in L^{2}(\mathbb{R}^{3})$ with $\nabla
\cdot u_{0}(x)=\nabla \cdot b_{0}(x)=0$ in last century, the problem on the
uniqueness and regularity of the weak solutions is one of the most
challenging problem of the mathematical community. Hence, many researchers
have developed different regularity criteria for the 3D MHD equations under
assumption of certain growth conditions on the velocity or on the magnetic
field (see, e.g., \cite{CW, CGG, CGZ, CMZ, DJZ, G1, HX, HW, JL,W1, W2, W3,
Z1, ZG1, ZG2} and the references therein).

Recent years, the problem of so-called regularity criteria via one
components\ was investigated for the MHD equations by some researchers (see
\cite{B, G2, G3, G4, JZ0, JZ1, JZ2, JZ3, JZ4, NGZ, Y16, Y14, Zh3} and the
references therein). In particular, in \cite{Y17}, Yamazaki established the
following regularity criterion by involving one velocity and one current
density component, which shows that a weak solution $(u,b)$ is smooth on a
time interval $(0,T]$ if
\begin{equation}
u_{3}\in L^{p}(0,T;L^{q}(\mathbb{R}^{3}))\text{ \ with \ }\frac{2}{p}+\frac{3%
}{q}\leq \frac{1}{3}+\frac{1}{2q},\text{ \ }\frac{15}{2}<q\leq \infty ,
\label{eq6}
\end{equation}%
and
\begin{equation*}
j_{3}\in L^{p^{\prime }}(0,T;L^{q^{\prime }}(\mathbb{R}^{3}))\text{ \ with \
}\frac{2}{p^{\prime }}+\frac{3}{q^{\prime }}\leq 2,\text{ \ }\frac{3}{2}%
<q^{\prime }\leq \infty ,
\end{equation*}%
where $j_{3}$ is the third component of the current density $j=\nabla \times
b=(j_{1},j_{2},j_{3})$. Later, Zhang \cite{Zh1} improved the regularity
criterion (\ref{eq6}) to the following conditions
\begin{equation}
u_{3}\in L^{p}(0,T;L^{q}(\mathbb{R}^{3}))\text{ \ with \ }\frac{2}{p}+\frac{3%
}{q}=\frac{4}{9}-\frac{1}{3q},\text{ \ }\frac{15}{2}\leq q\leq \infty ,
\label{eq7}
\end{equation}%
and
\begin{equation*}
j_{3}\in L^{p^{\prime }}(0,T;L^{q^{\prime }}(\mathbb{R}^{3}))\text{ \ with \
}\frac{2}{p^{\prime }}+\frac{3}{q^{\prime }}\leq 2,\text{ \ }\frac{3}{2}%
<q^{\prime }\leq \infty .
\end{equation*}%
Very recently, this result (\ref{eq7}) is further refined by Zhang \cite{Zh2}
to prove the regularity criterion as long as the following conditions
\begin{equation}
u_{3}\in L^{p}(0,T;L^{q}(\mathbb{R}^{3}))\text{ \ with \ \ }\frac{2}{p}+%
\frac{3}{q}=\frac{4}{9},\text{ \ }\frac{27}{4}\leq q\leq \infty ,
\label{eq08}
\end{equation}%
are satisfied.

Motivated by the papers \cite{Y17, Zh1, Zh2}, the purpose of the present
paper is to refine (\ref{eq08}) and to extend the above regularity criterion
to the Lorentz space $L^{\alpha ,\infty }$ which is larger than $L^{\alpha }$%
. More precisely, our main result now read as follows.

\begin{thm}
\label{th1}Suppose $T>0$, $(u_{0},b_{0})\in L^{2}(\mathbb{R}^{3})$ and $%
\nabla \cdot u_{0}=\nabla \cdot b_{0}=0$ in the sense of distributions.
Assume that $\left( u,b\right) $ is a weak solution of the 3D MHD equations (%
\ref{eq1.1}) on $(0,T)$.\ If $u_{3}$ and $j_{3}$ satisfy the following
growth conditions
\begin{equation}
\int_{0}^{T}(\left\Vert u_{3}(\tau )\right\Vert _{L^{\alpha ,\infty }}^{%
\frac{30\alpha }{7\alpha -45}}+\left\Vert j_{3}(\tau )\right\Vert _{L^{\beta
,\infty }}^{\frac{2\beta }{2\beta -3}})d\tau <\infty ,  \label{eq16}
\end{equation}%
where $\frac{45}{7}\leq \alpha \leq \infty $ and $\frac{3}{2}<\beta \leq
\infty $, then the weak solution $\left( u,b\right) $ is regular on $(0;T]$.
\end{thm}

\begin{re}
Theorem \ref{th1} extends the previous results on Navier-Stokes equations
due to the fact that the MHD equations with $b(x,t)=0$ reduces the
Navier-Stokes equations. According to the embedding relation $L^{\alpha
}\subseteq L^{\alpha ,\infty }$, it is easy to see that our result of
Theorem \ref{th1} is an improvement of the recent works by Yamazaki \cite%
{Y17} and Zhang  \cite{Zh1, Zh2}.
\end{re}

\section{Preliminaries}

Throughout this paper, we use the following usual notations. $L^{p}(\mathbb{R%
}^{3})$ denotes the Lebegue space associated with norm%
\begin{equation*}
\left\Vert f\right\Vert _{L^{p}}=\left\{
\begin{array}{c}
\left( \int_{\mathbb{R}^{3}}\left\vert f(x)\right\vert ^{p}dx\right) ^{\frac{%
1}{p}},\text{ \ \ for \ \ }1\leq p<\infty , \\
\underset{x\in \mathbb{R}^{3}}{ess\sup }\left\vert f(x)\right\vert ,\text{ \
\ for \ \ }p=\infty .%
\end{array}%
\right.
\end{equation*}%
$H^{k}(\mathbb{R}^{3})$ denotes the Hilbert space $\left\{ u\in L^{2}(%
\mathbb{R}^{3}):\left\Vert \nabla ^{k}u\right\Vert _{L^{2}}<\infty \right\} $%
. Let $(X,\mathcal{M},\mu )$ be a non-atomic measurable space. For a
complex- or real-valued $\mu -$measurable function $f(x)$ defined on $X$,
its distributional function is defined by%
\begin{equation*}
f_{\ast }(\sigma )=\mu \left\{ x\in X:f(x)>\sigma \right\} ,\text{ \ \ for \
}\sigma >0,
\end{equation*}%
which is non-increasing and continuous from the right. Furthermore, its
non-increasing rearrangement $f^{\ast }$ is defined by%
\begin{equation*}
f^{\ast }(t)=\inf \left\{ s>0:f_{\ast }(s)\leq t\right\} ,\text{ \ \ for \ }%
t>0,
\end{equation*}%
which is also non-increasing and continuous from the right and has the same
distributional function as $f(x)$.

The Lorentz space $L^{p,q}$ on $(X,\mathcal{M},\mu )$ is the collection of
all real- or complex-valued $\mu -$measurable functions $f(x)$ defined on $X$
such \ that $\left\Vert f\right\Vert _{L^{p,q}}<\infty $, where%
\begin{equation*}
\left\Vert f\right\Vert _{L^{p,q}}=\left\{
\begin{array}{c}
\left( \frac{q}{p}\int_{0}^{\infty }(t^{\frac{1}{p}}f^{\ast }(t))^{q}\frac{dt%
}{t}\right) ^{\frac{1}{q}},\text{ \ \ if \ \ }1\leq p<\infty ,\text{ \ }%
1<q<\infty  \\
\underset{t>0}{\sup }(t^{\frac{1}{p}}f^{\ast }(t)),\text{ \ \ if \ \ }1\leq
p<\infty ,\text{ \ }q=\infty .%
\end{array}%
\right.
\end{equation*}%
Moreover,%
\begin{equation*}
\left\Vert f\right\Vert _{L^{p,\infty }}=\underset{t>0}{\sup }(t^{\frac{1}{p}%
}f^{\ast }(t))=\underset{\sigma >0}{\sup }\sigma (f_{\ast }(\sigma ))^{\frac{%
1}{p}}
\end{equation*}%
for any $f\in L^{p,\infty }$. For details, we refer to \cite{BL} and \cite%
{Tri}.

The space definition implies the following continuous embeddings:%
\begin{equation*}
L^{p}(\mathbb{R}^{3})=L^{p,p}(\mathbb{R}^{3})\hookrightarrow L^{p,q}(\mathbb{%
R}^{3})\hookrightarrow L^{p,\infty }(\mathbb{R}^{3}),\text{ \ \ }1\leq p\leq
q<\infty .
\end{equation*}%
In order to prove Theorem \ref{th1}, we need the following H\"{o}lder
inequality in Lorentz spaces (see, e.g., O'Neil \cite{O} and \cite{Lor}).

\begin{lem}[\protect\cite{O}, Theorems 3.4 and 3.5]
\label{lem1}Let $f\in L^{p_{2},q_{2}}(\mathbb{R}^{3})$ and $g\in
L^{p_{3},q_{3}}(\mathbb{R}^{3})$ with $1\leq p_{2},p_{3}\leq \infty $, $%
1\leq q_{2},q_{3}\leq \infty $.\ Then $fg\in L^{p_{1},q_{1}}(\mathbb{R}^{3})$
with
\begin{equation*}
\frac{1}{p_{1}}=\frac{1}{p_{2}}+\frac{1}{p_{3}},\frac{1}{q_{1}}=\frac{1}{%
q_{2}}+\frac{1}{q_{3}}
\end{equation*}%
and the H\"{o}lder inequality of Lorentz spaces
\begin{equation*}
\left\Vert fg\right\Vert _{L^{p_{1},q_{1}}}\leq C\left\Vert f\right\Vert
_{L^{p_{2},q_{2}}}\left\Vert g\right\Vert _{L^{p_{3},q_{3}}},
\end{equation*}%
holds true for a positive constant $C$.
\end{lem}

We also recall Gagliardo-Nirenberg's inequality in Lorentz spaces which
plays an important role in the proofs of Theorem \ref{th1}.

\begin{lem}
\label{lem2}Let $f\in L^{p,q}(\mathbb{R}^{3})$ with $1\leq
p,q,p_{4},q_{4},p_{5},q_{5}\leq \infty $. Then the Gagliardo-Nirenberg
inequality of Lorentz spaces%
\begin{equation*}
\left\Vert f\right\Vert _{L^{p,q}}\leq C\left\Vert f\right\Vert
_{L^{p_{4},q_{4}}}^{\theta }\left\Vert f\right\Vert
_{L^{p_{5},q_{5}}}^{1-\theta }
\end{equation*}%
holds for a positive constant $C$ and%
\begin{equation*}
\frac{1}{p}=\frac{\theta }{p_{4}}+\frac{1-\theta }{p_{5}},\text{ \ \ }\frac{1%
}{q}=\frac{\theta }{q_{4}}+\frac{1-\theta }{q_{5}},\text{ \ }\theta \in
(0,1).
\end{equation*}
\end{lem}

\section{Proof of main result.}

In this section, under the assumptions of the Theorem \ref{th1}, we prove
our main result. Before proving our result, we recall the following
muliplicative Sobolev imbedding inequality in the whole space $\mathbb{R}%
^{3} $ (see, for example \cite{CT}) :
\begin{equation}
\left\Vert f\right\Vert _{L^{6}}\leq C\left\Vert \nabla _{h}f\right\Vert
_{L^{2}}^{\frac{2}{3}}\left\Vert \partial _{3}f\right\Vert _{L^{2}}^{\frac{1%
}{3}},  \label{eq8}
\end{equation}%
where $\nabla _{h}=(\partial _{x_{1}},\partial _{x_{2}})$ is the horizontal
gradient operator. We are now give the proof of our main theorem.

\begin{pf}
To prove our result, it suffices to show that for any fixed $T>T^{\ast }$,
there holds%
\begin{equation*}
\underset{0\leq t\leq T^{\ast }}{\sup }\left\Vert \nabla u(t)\right\Vert
_{L^{2}}^{2}+\left\Vert \nabla b(t)\right\Vert _{L^{2}}^{2}\leq C_{T},
\end{equation*}%
where $T^{\ast }$, which denotes the maximal existence time of a strong
solution and $C_{T}$ is an absolute constant which only depends on $T,u_{0}$
and $b_{0}$.

The method of our proof is the standard energy estimates as in \cite{Y17}.
We will based on two major parts. The first one establishes the bounds of $%
(\left\Vert \nabla _{h}u\right\Vert _{L^{2}}^{2}+\left\Vert \nabla
_{h}b\right\Vert _{L^{2}}^{2})$, while the second gives the bounds of the $%
H^{1}-$norm of velocity $u$ and magnetic field $b$ in terms of the results
of part one.

For this purpose, we multiply the first and second equations of (\ref{eq1.1}%
) by $-\Delta _{h}u$ and $-\Delta _{h}b$, respectively, and integrate them
over $\mathbb{R}^{3}$ with respect to the spatial variable. Then,
integration by parts gives the following identity:
\begin{eqnarray*}
&&\frac{1}{2}\frac{d}{dt}(\left\Vert \nabla _{h}u\right\Vert
_{L^{2}}^{2}+\left\Vert \nabla _{h}b\right\Vert _{L^{2}}^{2})+\left\Vert
\nabla \nabla _{h}u\right\Vert _{L^{2}}^{2}+\left\Vert \nabla \nabla
_{h}b\right\Vert _{L^{2}}^{2} \\
&=&\int_{\mathbb{R}^{3}}(u\cdot \nabla )u\cdot \Delta _{h}udx-\int_{\mathbb{R%
}^{3}}(b\cdot \nabla )b\cdot \Delta _{h}udx \\
&&+\int_{\mathbb{R}^{3}}(u\cdot \nabla )b\cdot \Delta _{h}bdx-\int_{\mathbb{R%
}^{3}}(b\cdot \nabla )u\cdot \Delta _{h}bdx \\
&=&\text{RHS},
\end{eqnarray*}%
where $\Delta _{h}=\partial _{x_{1}}^{2}+\partial _{x_{2}}^{2}$ is the
horizontal Laplacian. For simplicity of exposition, we denote%
\begin{eqnarray*}
\mathcal{L}^{2}(t) &=&\underset{\tau \in \lbrack \Gamma ,t]}{\sup }%
(\left\Vert \nabla _{h}u(\tau )\right\Vert _{L^{2}}^{2}+\left\Vert \nabla
_{h}b(\tau )\right\Vert _{L^{2}}^{2})+\int_{\Gamma }^{t}(\left\Vert \nabla
\nabla _{h}u(\tau )\right\Vert _{L^{2}}^{2}+\left\Vert \nabla \nabla
_{h}b(\tau )\right\Vert _{L^{2}}^{2})d\tau , \\
\mathcal{J}^{2}(t) &=&\underset{\tau \in \lbrack \Gamma ,t]}{\sup }%
(\left\Vert \nabla u(\tau )\right\Vert _{L^{2}}^{2}+\left\Vert \nabla b(\tau
)\right\Vert _{L^{2}}^{2})+\int_{\Gamma }^{t}(\left\Vert \Delta u(\tau
)\right\Vert _{L^{2}}^{2}+\left\Vert \Delta b(\tau )\right\Vert
_{L^{2}}^{2})d\tau ,
\end{eqnarray*}%
for $t\in \lbrack \Gamma ,T^{\ast })$. We choose $\epsilon ,\eta >0$ to be
precisely determined subsequently and then select $\Gamma <T^{\ast }$
sufficiently close to $T^{\ast }$ such that for all $\Gamma \leq t<T^{\ast }$%
,%
\begin{equation}
\int_{\Gamma }^{t}(\left\Vert \nabla u(\tau )\right\Vert
_{L^{2}}^{2}+\left\Vert \nabla b(\tau )\right\Vert _{L^{2}}^{2})d\tau \leq
\epsilon \ll 1\text{ \ \ \ and  \ \ }\int_{\Gamma }^{t}\left\Vert j_{3}(\tau
)\right\Vert _{L^{\beta }}^{\frac{2\beta }{2\beta -3}}d\tau \leq \eta \ll 1.
\label{eq19}
\end{equation}%
Applying the divergence-free condition, $\nabla \cdot u=\nabla \cdot b=0$,
we find that RHS can be estimated as
\begin{eqnarray}
\text{RHS} &\leq &\int_{\mathbb{R}^{3}}\left\vert u_{3}\right\vert
\left\vert \nabla u\right\vert \left\vert \nabla \nabla _{h}u\right\vert
dx+\int_{\mathbb{R}^{3}}\left\vert u_{3}\right\vert \left\vert \nabla
b\right\vert \left\vert \nabla \nabla _{h}b\right\vert dx+\int_{\mathbb{R}%
^{3}}\left\vert b_{3}\right\vert \left\vert \nabla u\right\vert \left\vert
\nabla \nabla _{h}b\right\vert dx  \notag \\
&&+\int_{\mathbb{R}^{3}}\left\vert b_{3}\right\vert \left\vert \nabla
b\right\vert \left\vert \nabla \nabla _{h}u\right\vert dx+\int_{\mathbb{R}%
^{3}}\left\vert \nabla _{h}u\right\vert \left\vert \nabla _{h}b\right\vert
\left\vert j_{3}\right\vert dx  \notag \\
&=&L_{1}+L_{2}+L_{3}+L_{4}+L_{5},  \label{eq3}
\end{eqnarray}%
where the last inequality was proved in \cite{Y17} (see, Proposition 3.1 in
\cite{Y17} for details).

With the use of the Lemma \ref{lem1}, (\ref{eq8}), and the Young inequality,
we derive the estimate of the first term $L_{1}$ of (\ref{eq3}) as follows :%
\begin{eqnarray*}
L_{1} &\leq &C\left\Vert u_{3}\right\Vert _{L^{\alpha ,\infty }}\left\Vert
\nabla u\right\Vert _{L^{\frac{2\alpha }{\alpha -2},2}}\left\Vert \nabla
\nabla _{h}u\right\Vert _{L^{2}} \\
&\leq &C\left\Vert u_{3}\right\Vert _{L^{\alpha ,\infty }}\left\Vert \nabla
u\right\Vert _{L^{2}}^{1-\frac{3}{\alpha }}\left\Vert \nabla u\right\Vert
_{L^{6}}^{\frac{3}{\alpha }}\left\Vert \nabla \nabla _{h}u\right\Vert
_{L^{2}} \\
&\leq &C\left\Vert u_{3}\right\Vert _{L^{\alpha ,\infty }}\left\Vert \nabla
u\right\Vert _{L^{2}}^{1-\frac{3}{\alpha }}\left\Vert \Delta u\right\Vert
_{L^{2}}^{\frac{1}{\alpha }}\left\Vert \nabla \nabla _{h}u\right\Vert
_{L^{2}}^{1+\frac{2}{\alpha }} \\
&\leq &C\left\Vert u_{3}\right\Vert _{L^{\alpha ,\infty }}^{\frac{2\alpha }{%
\alpha -2}}\left\Vert \nabla u\right\Vert _{L^{2}}^{2-\frac{2}{\alpha -2}%
}\left\Vert \Delta u\right\Vert _{L^{2}}^{\frac{2}{\alpha -2}}+\frac{1}{8}%
\left\Vert \nabla \nabla _{h}u\right\Vert _{L^{2}}^{2},
\end{eqnarray*}%
where we have used the following Gagliardo-Nirenberg inequality in Lorentz
spaces :
\begin{equation*}
\left\Vert \nabla u\right\Vert _{L^{\frac{2\alpha }{\alpha -2},2}}\leq
C\left\Vert \nabla u\right\Vert _{L^{2}}^{1-\frac{3}{\alpha }}\left\Vert
\nabla u\right\Vert _{L^{6}}^{\frac{3}{\alpha }}.
\end{equation*}%
Similarly, employing the H\"{o}lder inequality and the Gagliardo-Nirenberg
inequality give that for $L_{1}$,
\begin{equation*}
L_{2}\leq C\left\Vert u_{3}\right\Vert _{L^{\alpha ,\infty }}^{\frac{2\alpha
}{\alpha -2}}\left\Vert \nabla b\right\Vert _{L^{2}}^{2-\frac{2}{\alpha -2}%
}\left\Vert \Delta b\right\Vert _{L^{2}}^{\frac{2}{\alpha -2}}+\frac{1}{8}%
\left\Vert \nabla \nabla _{h}b\right\Vert _{L^{2}}^{2}.
\end{equation*}%
We now estimate $L_{3}$,
\begin{eqnarray*}
L_{3} &\leq &\left\Vert b_{3}\right\Vert _{L^{10}}\left\Vert \nabla
u\right\Vert _{L^{\frac{5}{2}}}\left\Vert \nabla \nabla _{h}b\right\Vert
_{L^{2}}\leq C\left\Vert b_{3}\right\Vert _{L^{10}}\left\Vert \nabla
u\right\Vert _{L^{2}}^{\frac{7}{10}}\left\Vert \nabla u\right\Vert _{L^{6}}^{%
\frac{3}{10}}\left\Vert \nabla \nabla _{h}b\right\Vert _{L^{2}} \\
&\leq &C\left\Vert b_{3}\right\Vert _{L^{\frac{10}{3}}}\left\Vert \nabla
u\right\Vert _{L^{2}}^{\frac{7}{10}}\left\Vert \nabla \nabla
_{h}u\right\Vert _{L^{2}}^{\frac{1}{5}}\left\Vert \Delta u\right\Vert
_{L^{2}}^{\frac{1}{10}}\left\Vert \nabla \nabla _{h}b\right\Vert _{L^{2}} \\
&\leq &C\left\Vert b_{3}\right\Vert _{L^{10}}^{\frac{5}{2}}\left\Vert \nabla
u\right\Vert _{L^{2}}^{\frac{7}{4}}\left\Vert \Delta u\right\Vert _{L^{2}}^{%
\frac{1}{4}}+\frac{1}{8}(\left\Vert \nabla \nabla _{h}b\right\Vert
_{L^{2}}^{2}+\left\Vert \nabla \nabla _{h}u\right\Vert _{L^{2}}^{2}),
\end{eqnarray*}%
where we have used the fact$\left\Vert \nabla u\right\Vert _{L^{\frac{5}{2}%
}}\leq C\left\Vert \nabla u\right\Vert _{L^{2}}^{\frac{7}{10}}\left\Vert
\nabla u\right\Vert _{L^{6}}^{\frac{3}{10}}$.

Likewise,
\begin{equation*}
L_{4}\leq C\left\Vert b_{3}\right\Vert _{L^{10}}^{\frac{5}{2}}\left\Vert
\nabla b\right\Vert _{L^{2}}^{\frac{7}{4}}\left\Vert \Delta b\right\Vert
_{L^{2}}^{\frac{1}{4}}+\frac{1}{8}(\left\Vert \nabla \nabla _{h}b\right\Vert
_{L^{2}}^{2}+\left\Vert \nabla \nabla _{h}u\right\Vert _{L^{2}}^{2}).
\end{equation*}%
For $L_{5}$, by applying the H\"{o}lder inequality, the Gagliardo-Nirenberg
inequality and the Young inequality, one shows that
\begin{eqnarray*}
L_{5} &=&\int_{\mathbb{R}^{3}}\left\vert \nabla _{h}u\right\vert \left\vert
\nabla _{h}b\right\vert \left\vert j_{3}\right\vert dx\leq \frac{1}{2}\int_{%
\mathbb{R}^{3}}\left( \left\vert \nabla _{h}u\right\vert ^{2}+\left\vert
\nabla _{h}b\right\vert ^{2}\right) \left\vert j_{3}\right\vert dx \\
&\leq &C\left\Vert j_{3}\right\Vert _{L^{\beta ,\infty }}(\left\Vert \nabla
_{h}u\right\Vert _{L^{\frac{2\beta }{\beta -2},2}}\left\Vert \nabla
_{h}u\right\Vert _{L^{2}}+\left\Vert \nabla _{h}u\right\Vert _{L^{\frac{%
2\beta }{\beta -2},2}}\left\Vert \nabla _{h}b\right\Vert _{L^{2}}) \\
&\leq &C\left\Vert j_{3}\right\Vert _{L^{\beta ,\infty }}(\left\Vert \nabla
_{h}u\right\Vert _{L^{2}}^{2-\frac{3}{\beta }}\left\Vert \nabla \nabla
_{h}u\right\Vert _{L^{2}}^{\frac{3}{s}}+\left\Vert \nabla _{h}b\right\Vert
_{L^{2}}^{2-\frac{3}{\beta }}\left\Vert \nabla \nabla _{h}b\right\Vert
_{L^{2}}^{\frac{3}{s}}) \\
&\leq &C\left\Vert j_{3}\right\Vert _{L^{\beta ,\infty }}^{\frac{2\beta }{%
2\beta -3}}(\left\Vert \nabla _{h}u\right\Vert _{L^{2}}^{2}+\left\Vert
\nabla _{h}b\right\Vert _{L^{2}}^{2})+\frac{1}{8}(\left\Vert \nabla \nabla
_{h}b\right\Vert _{L^{2}}^{2}+\left\Vert \nabla \nabla _{h}u\right\Vert
_{L^{2}}^{2}).
\end{eqnarray*}%
Inserting all the estimates into (\ref{eq3}), Gronwall's type argument using
\begin{equation*}
1\leq \underset{\lambda \in \lbrack \Gamma ,\tau ]}{\sup }\exp \left(
c\int_{\lambda }^{\tau }\left\Vert j_{3}(\varphi )\right\Vert _{L^{\beta
,\infty }}^{\frac{2\beta }{2\beta -3}}d\varphi \right) \lesssim \exp \left(
c\int_{0}^{T^{\ast }}\left\Vert j_{3}(\varphi )\right\Vert _{L^{\beta
,\infty }}^{\frac{2\beta }{2\beta -3}}d\varphi \right) \lesssim 1,
\end{equation*}%
due to (\ref{eq16}) leads to, for every $\tau \in \lbrack \Gamma ,t]$%
\begin{eqnarray}
\mathcal{L}^{2}(t) &\leq &C+C\int_{\Gamma }^{t}\left\Vert u_{3}\right\Vert
_{L^{\alpha ,\infty }}^{\frac{2\alpha }{\alpha -2}}(\left\Vert \nabla
u\right\Vert _{L^{2}}^{2-\frac{2}{\alpha -2}}\left\Vert \Delta u\right\Vert
_{L^{2}}^{\frac{2}{\alpha -2}}+\left\Vert \nabla b\right\Vert _{L^{2}}^{2-%
\frac{2}{\alpha -2}}\left\Vert \Delta b\right\Vert _{L^{2}}^{\frac{2}{\alpha
-2}})d\tau   \notag \\
&&+C\int_{\Gamma }^{t}\left\Vert b_{3}\right\Vert _{L^{10}}^{\frac{5}{2}%
}(\left\Vert \nabla u\right\Vert _{L^{2}}^{\frac{7}{4}}\left\Vert \Delta
u\right\Vert _{L^{2}}^{\frac{1}{4}}+\left\Vert \nabla b\right\Vert _{L^{2}}^{%
\frac{7}{4}}\left\Vert \Delta b\right\Vert _{L^{2}}^{\frac{1}{4}})d\tau
\notag \\
&&+C\int_{\Gamma }^{t}\left\Vert j_{3}\right\Vert _{L^{\beta ,\infty }}^{%
\frac{2\beta }{2\beta -3}}(\left\Vert \nabla _{h}u\right\Vert
_{L^{2}}^{2}+\left\Vert \nabla _{h}b\right\Vert _{L^{2}}^{2})d\tau   \notag
\\
&=&C+\mathcal{I}_{1}(t)+\mathcal{I}_{2}(t)+\mathcal{I}_{3}(t).  \label{eq4}
\end{eqnarray}%
Next, we analyze the right-hand side of (\ref{eq4}) one by one. First, due
to (\ref{eq19}) and the definition of $\mathcal{J}^{2}$, we have
\begin{eqnarray*}
\mathcal{I}_{1}(t) &\leq &C\left( \underset{\tau \in \lbrack \Gamma ,t]}{%
\sup }\left\Vert \nabla u(\tau )\right\Vert _{L^{2}}^{\frac{3}{2}-\frac{2}{%
\alpha -2}}\right) \int_{\Gamma }^{t}\left\Vert u_{3}(\tau )\right\Vert
_{L^{\alpha ,\infty }}^{\frac{2\alpha }{\alpha -2}}\left\Vert \nabla u(\tau
)\right\Vert _{L^{2}}^{\frac{1}{2}}\left\Vert \Delta u(\tau )\right\Vert
_{L^{2}}^{\frac{2}{\alpha -2}}d\tau  \\
&&+C\left( \underset{\tau \in \lbrack \Gamma ,t]}{\sup }\left\Vert \nabla
b(\tau )\right\Vert _{L^{2}}^{\frac{3}{2}-\frac{2}{\alpha -2}}\right)
\int_{\Gamma }^{t}\left\Vert u_{3}(\tau )\right\Vert _{L^{\alpha ,\infty }}^{%
\frac{2\alpha }{\alpha -2}}\left\Vert \nabla b(\tau )\right\Vert _{L^{2}}^{%
\frac{1}{2}}\left\Vert \Delta b(\tau )\right\Vert _{L^{2}}^{\frac{2}{\alpha
-2}}d\tau  \\
&\leq &C\mathcal{J}^{\frac{3}{2}-\frac{2}{\alpha -2}}(t)\left( \int_{\Gamma
}^{t}\left\Vert u_{3}(\tau )\right\Vert _{L^{\alpha ,\infty }}^{\frac{%
8\alpha }{3\alpha -10}}d\tau \right) ^{\frac{3}{4}-\frac{1}{\alpha -2}%
}\left( \int_{\Gamma }^{t}\left\Vert \nabla u(\tau )\right\Vert
_{L^{2}}^{2}d\tau \right) ^{\frac{1}{4}}\left( \int_{\Gamma }^{t}\left\Vert
\Delta u(\tau )\right\Vert _{L^{2}}^{2}d\tau \right) ^{\frac{1}{\alpha -2}}
\\
&&+C\mathcal{J}^{\frac{3}{2}-\frac{2}{\alpha -2}}(t)\left( \int_{\Gamma
}^{t}\left\Vert u_{3}(\tau )\right\Vert _{L^{\alpha ,\infty }}^{\frac{%
8\alpha }{3\alpha -10}}d\tau \right) ^{\frac{3}{4}-\frac{1}{\alpha -2}%
}\left( \int_{\Gamma }^{t}\left\Vert \nabla b(\tau )\right\Vert
_{L^{2}}^{2}d\tau \right) ^{\frac{1}{4}}\left( \int_{\Gamma }^{t}\left\Vert
\Delta b(\tau )\right\Vert _{L^{2}}^{2}d\tau \right) ^{\frac{1}{\alpha -2}}
\\
&\leq &C\mathcal{J}^{\frac{3}{2}-\frac{2}{\alpha -2}}(t)\left( \int_{\Gamma
}^{t}\left\Vert u_{3}(\tau )\right\Vert _{L^{\alpha ,\infty }}^{\frac{%
8\alpha }{3\alpha -10}}d\tau \right) ^{\frac{3}{4}-\frac{1}{\alpha -2}%
}\epsilon ^{\frac{1}{4}}\mathcal{J}^{\frac{2}{\alpha -2}}(t) \\
&=&C\epsilon ^{\frac{1}{4}}\mathcal{J}^{\frac{3}{2}}(t)\left( \int_{\Gamma
}^{t}\left\Vert u_{3}(\tau )\right\Vert _{L^{\alpha ,\infty }}^{\frac{%
8\alpha }{3\alpha -10}}d\tau \right) ^{\frac{3}{4}-\frac{1}{\alpha -2}}.
\end{eqnarray*}%
Now, we estimate the term $\mathcal{I}_{2}(t)$ as
\begin{eqnarray*}
\mathcal{I}_{2}(t) &\leq &C\left( \underset{\tau \in \lbrack \Gamma ,t]}{%
\sup }\left\Vert b_{3}(\tau )\right\Vert _{L^{10}}^{\frac{5}{2}}\right)
\int_{\Gamma }^{t}\left\Vert \nabla u(\tau )\right\Vert _{L^{2}}^{\frac{7}{4}%
}\left\Vert \Delta u(\tau )\right\Vert _{L^{2}}^{\frac{1}{4}}d\tau  \\
&&+\left( \underset{\tau \in \lbrack \Gamma ,t]}{\sup }\left\Vert b_{3}(\tau
)\right\Vert _{L^{10}}^{\frac{5}{2}}\right) \int_{\Gamma }^{t}\left\Vert
\nabla b(\tau )\right\Vert _{L^{2}}^{\frac{7}{4}}\left\Vert \Delta b(\tau
)\right\Vert _{L^{2}}^{\frac{1}{4}}d\tau  \\
&\leq &\left( \underset{\tau \in \lbrack \Gamma ,t]}{\sup }\left\Vert
b_{3}(\tau )\right\Vert _{L^{10}}^{\frac{5}{2}}\right) \left( \int_{\Gamma
}^{t}\left\Vert \nabla u(\tau )\right\Vert _{L^{2}}^{2}d\tau \right) ^{\frac{%
7}{8}}\left( \int_{\Gamma }^{t}\left\Vert \Delta u(\tau )\right\Vert
_{L^{2}}^{2}d\tau \right) ^{\frac{1}{8}} \\
&&+\left( \underset{\tau \in \lbrack \Gamma ,t]}{\sup }\left\Vert b_{3}(\tau
)\right\Vert _{L^{10}}^{\frac{5}{2}}\right) \left( \int_{\Gamma
}^{t}\left\Vert \nabla b(\tau )\right\Vert _{L^{2}}^{2}d\tau \right) ^{\frac{%
7}{8}}\left( \int_{\Gamma }^{t}\left\Vert \Delta b(\tau )\right\Vert
_{L^{2}}^{2}d\tau \right) ^{\frac{1}{8}} \\
&\leq &C\left( \underset{\tau \in \lbrack \Gamma ,t]}{\sup }\left\Vert
b_{3}(\tau )\right\Vert _{L^{10}}^{\frac{5}{2}}\right) \epsilon ^{\frac{7}{8}%
}\mathcal{J}^{\frac{1}{4}}(t).
\end{eqnarray*}%
For $\mathcal{I}_{3}(t)$, applying H\"{o}lder's and Young's inequalities, we
get
\begin{eqnarray*}
\mathcal{I}_{3}(t) &\leq &C\underset{\tau \in \lbrack \Gamma ,t]}{\sup }%
(\left\Vert \nabla _{h}u(\tau )\right\Vert _{L^{2}}^{2}+\left\Vert \nabla
_{h}b(\tau )\right\Vert _{L^{2}}^{2})\int_{\Gamma }^{t}\left\Vert j_{3}(\tau
)\right\Vert _{L^{\beta ,\infty }}^{\frac{2\beta }{2\beta -3}}d\tau  \\
&\leq &C\eta \mathcal{L}^{2}(t).
\end{eqnarray*}%
Therefore, combining the estimates of $\mathcal{I}_{1}(t),\mathcal{I}_{2}(t)$
and $\mathcal{I}_{3}(t)$\ together with (\ref{eq4}) and taking $\eta $ small
enough, it is easy to see that for all $\Gamma \leq t<T^{\ast }:$%
\begin{equation}
\mathcal{L}^{2}(t)\leq C+C\epsilon ^{\frac{1}{4}}\mathcal{J}^{\frac{3}{2}%
}(t)\left( \int_{\Gamma }^{t}\left\Vert u_{3}(\tau )\right\Vert _{L^{\alpha
,\infty }}^{\frac{8\alpha }{3\alpha -10}}d\tau \right) ^{\frac{3\alpha -10}{%
4(\alpha -2)}}+C\left( \underset{\tau \in \lbrack \Gamma ,t]}{\sup }%
\left\Vert b_{3}(\tau )\right\Vert _{L^{10}}^{\frac{5}{2}}\right) \epsilon ^{%
\frac{7}{8}}\mathcal{J}^{\frac{1}{4}}(t)  \label{eq20}
\end{equation}

Now, we will establish the bounds of $L^{10}$-norm of the magnetic field $%
b_{3}$. In order to do it, we recall the third equation of the magnetic
field:%
\begin{equation*}
\partial _{t}b_{3}-\Delta b_{3}+(u\cdot \nabla )b_{3}=(b\cdot \nabla )u_{3},
\end{equation*}%
and multiply this equation by $\left\vert b_{3}\right\vert ^{8}b_{3}$,
integrating by parts, using incompressibility conditions to obtain%
\begin{eqnarray}
\frac{1}{10}\frac{d}{dt}\int_{\mathbb{R}^{3}}\left\vert b_{3}\right\vert
^{10}dx+\frac{9}{25}\int_{\mathbb{R}^{3}}\left\vert \nabla
(b_{3}^{5})\right\vert ^{2}dx &=&\int_{\mathbb{R}^{3}}(b\cdot \nabla
u_{3})(\left\vert b_{3}\right\vert ^{8}b_{3})dx  \notag \\
&=&-9\int_{\mathbb{R}^{3}}b\cdot \left\vert b_{3}\right\vert ^{4}(\left\vert
b_{3}\right\vert ^{4}\nabla b_{3})u_{3}dx  \notag \\
&\leq &\frac{9}{5}\int_{\mathbb{R}^{3}}\left\vert b\right\vert (\left\vert
b_{3}\right\vert ^{5})^{\frac{4}{5}}\left\vert u_{3}\right\vert \left\vert
\nabla (b_{3}^{5})\right\vert dx=I.  \label{eq22}
\end{eqnarray}%
Using the H\"{o}lder, Young inequalities and interpolation, the estimates of
$I$ is given by
\begin{eqnarray*}
I &\leq &\frac{9}{5}\left\Vert b\right\Vert _{L^{6}}\left\Vert u_{3}\left(
\left\vert b_{3}\right\vert ^{5}\right) ^{\frac{4}{5}}\right\Vert
_{L^{3}}\left\Vert \nabla (b_{3}^{5})\right\Vert _{L^{2}} \\
&\leq &C\left\Vert \nabla _{h}b\right\Vert _{L^{2}}^{\frac{2}{3}}\left\Vert
\nabla b\right\Vert _{L^{2}}^{\frac{1}{3}}\left\Vert u_{3}\right\Vert
_{L^{\alpha ,\infty }}\left\Vert \left( \left\vert b_{3}\right\vert
^{5}\right) ^{\frac{4}{5}}\right\Vert _{L^{\frac{3\alpha }{\alpha -3}%
,3}}\left\Vert \nabla (b_{3}^{5})\right\Vert _{L^{2}} \\
&\leq &C\left\Vert \nabla _{h}b\right\Vert _{L^{2}}^{\frac{2}{3}}\left\Vert
\nabla b\right\Vert _{L^{2}}^{\frac{1}{3}}\left\Vert u_{3}\right\Vert
_{L^{\alpha ,\infty }}\left\Vert \left( \left\vert b_{3}\right\vert
^{5}\right) ^{\frac{4}{5}}\right\Vert _{L^{\frac{3\alpha }{\alpha -3}%
,2}}\left\Vert \nabla (b_{3}^{5})\right\Vert _{L^{2}} \\
&\leq &C\left\Vert \nabla _{h}b\right\Vert _{L^{2}}^{\frac{2}{3}}\left\Vert
\nabla b\right\Vert _{L^{2}}^{\frac{1}{3}}\left\Vert u_{3}\right\Vert
_{L^{\alpha ,\infty }}\left( \left\Vert b_{3}^{5}\right\Vert _{L^{2}}^{\frac{%
3(\alpha -5)}{4\alpha }}\left\Vert \nabla (b_{3}^{5})\right\Vert _{L^{2}}^{%
\frac{15+\alpha }{4\alpha }}\right) ^{\frac{4}{5}}\left\Vert \nabla
(b_{3}^{5})\right\Vert _{L^{2}} \\
&\leq &C\left\Vert \nabla _{h}b\right\Vert _{L^{2}}^{\frac{20\alpha }{%
3(4\alpha -15)}}\left\Vert \nabla b\right\Vert _{L^{2}}^{\frac{10\alpha }{%
3(4\alpha -15)}}\left\Vert u_{3}\right\Vert _{L^{\alpha }}^{\frac{10\alpha }{%
4\alpha -15}}\left\Vert b_{3}^{5}\right\Vert _{L^{2}}^{\frac{6(\alpha -5)}{%
4\alpha -15}}+\frac{9}{5}\left\Vert \nabla (b_{3}^{5})\right\Vert
_{L^{2}}^{2}
\end{eqnarray*}%
Putting $I$ in (\ref{eq22}), we get
\begin{equation*}
\frac{d}{dt}\int_{\mathbb{R}^{3}}\left\vert b_{3}\right\vert ^{10}dx\leq
C\left\Vert \nabla _{h}b\right\Vert _{L^{2}}^{\frac{20\alpha }{3(4\alpha -15)%
}}\left\Vert \nabla b\right\Vert _{L^{2}}^{\frac{10\alpha }{3(4\alpha -15)}%
}\left\Vert u_{3}\right\Vert _{L^{\alpha }}^{\frac{10\alpha }{4\alpha -15}%
}\left\Vert b_{3}\right\Vert _{L^{10}}^{\frac{30(\alpha -5)}{4\alpha -15}}.
\end{equation*}%
Dividing by $\left\Vert b_{3}\right\Vert _{L^{10}}^{\frac{30(\alpha -5)}{%
4\alpha -15}}$, we arrive at%
\begin{equation*}
\frac{d}{dt}\left\Vert b_{3}\right\Vert _{L^{10}}^{\frac{10\alpha }{4\alpha
-15}}\leq C\left\Vert \nabla _{h}b\right\Vert _{L^{2}}^{\frac{20\alpha }{%
3(7\alpha -15)}}\left\Vert \nabla b\right\Vert _{L^{2}}^{\frac{10\alpha }{%
3(4\alpha -15)}}\left\Vert u_{3}\right\Vert _{L^{\alpha }}^{\frac{10\alpha }{%
4\alpha -15}}.
\end{equation*}%
Integrating over interval $[\Gamma ,\tau )$, it follows that%
\begin{equation}
\left\Vert b_{3}(\tau )\right\Vert _{L^{10}}\leq \left[ \left\Vert
b_{3}(\Gamma )\right\Vert _{L^{10}}^{\frac{10\alpha }{4\alpha -15}%
}+C\int_{\Gamma }^{\tau }\left\Vert \nabla _{h}b(\lambda )\right\Vert
_{L^{2}}^{\frac{20\alpha }{3(4\alpha -15)}}\left\Vert \nabla b(\lambda
)\right\Vert _{L^{2}}^{\frac{10\alpha }{3(4\alpha -15)}}\left\Vert
u_{3}(\lambda )\right\Vert _{L^{\alpha }}^{\frac{10\alpha }{4\alpha -15}%
}d\lambda \right] ^{\frac{4\alpha -15}{10\alpha }},  \label{eq21}
\end{equation}%
for all $\tau \in \lbrack \Gamma ,t)$. It follows from (\ref{eq21}) and (\ref%
{eq20}) that%
\begin{eqnarray*}
\mathcal{L}^{2}(t) &\leq &C+C\epsilon ^{\frac{1}{4}}\mathcal{J}^{\frac{3}{2}%
}(t)\left( \int_{\Gamma }^{t}\left\Vert u_{3}(\tau )\right\Vert _{L^{\alpha
}}^{\frac{8\alpha }{3\alpha -10}}d\tau \right) ^{\frac{3\alpha -10}{4(\alpha
-2)}} \\
&&+C\epsilon ^{\frac{7}{8}}\mathcal{J}^{\frac{1}{4}}(t)\underset{\tau \in
\lbrack \Gamma ,t]}{\sup }\left[ \left\Vert b_{3}(\Gamma )\right\Vert
_{L^{10}}^{\frac{10\alpha }{7\alpha -15}}+C\int_{\Gamma }^{\tau }\left\Vert
\nabla _{h}b(\lambda )\right\Vert _{L^{2}}^{\frac{20\alpha }{3(4\alpha -15)}%
}\left\Vert \nabla b(\lambda )\right\Vert _{L^{2}}^{\frac{10\alpha }{%
3(4\alpha -15)}}\left\Vert u_{3}(\lambda )\right\Vert _{L^{\alpha }}^{\frac{%
10\alpha }{4\alpha -15}}d\lambda \right] ^{\frac{4\alpha -15}{4\alpha }} \\
&\leq &C+C\epsilon ^{\frac{1}{4}}\mathcal{J}^{\frac{3}{2}}(t)\left(
\int_{\Gamma }^{t}\left\Vert u_{3}(\tau )\right\Vert _{L^{\alpha }}^{\frac{%
8\alpha }{3\alpha -10}}d\tau \right) ^{\frac{3\alpha -10}{4(\alpha -2)}} \\
&&+C\epsilon ^{\frac{7}{8}}\mathcal{J}^{\frac{1}{4}}(t)\left[ \left\Vert
b_{3}(\Gamma )\right\Vert _{L^{10}}^{\frac{10\alpha }{4\alpha -15}%
}+C\int_{\Gamma }^{t}\left\Vert \nabla _{h}b(\tau )\right\Vert _{L^{2}}^{%
\frac{20\alpha }{3(4\alpha -15)}}\left\Vert \nabla b(\tau )\right\Vert
_{L^{2}}^{\frac{10\alpha }{3(4\alpha -15)}}\left\Vert u_{3}(\tau
)\right\Vert _{L^{\alpha }}^{\frac{10\alpha }{4\alpha -15}}d\tau \right] ^{%
\frac{4\alpha -15}{4\alpha }} \\
&\leq &C+C\epsilon ^{\frac{1}{4}}\mathcal{J}^{\frac{3}{2}}(t)\left(
\int_{\Gamma }^{t}\left\Vert u_{3}(\tau )\right\Vert _{L^{\alpha }}^{\frac{%
8\alpha }{3\alpha -10}}d\tau \right) ^{\frac{3\alpha -10}{4(\alpha -2)}%
}+C\left\Vert b_{3}(\Gamma )\right\Vert _{L^{10}}^{\frac{5}{2}}\epsilon ^{%
\frac{7}{8}}\mathcal{J}^{\frac{1}{4}}(t) \\
&&+C\epsilon ^{\frac{7}{8}}\mathcal{J}^{\frac{1}{4}}(t)\underset{\tau \in
\lbrack \Gamma ,t]}{\sup }\left\Vert \nabla _{h}b(\tau )\right\Vert
_{L^{2}}^{\frac{5}{3}}\left[ \int_{\Gamma }^{t}\left\Vert \nabla b(\tau
)\right\Vert _{L^{2}}^{\frac{10\alpha }{3(4\alpha -15)}}\left\Vert
u_{3}(\tau )\right\Vert _{L^{\alpha }}^{\frac{10\alpha }{4\alpha -15}}d\tau %
\right] ^{\frac{4\alpha -15}{4\alpha }} \\
&\leq &C+C\epsilon ^{\frac{1}{4}}\mathcal{J}^{\frac{3}{2}}(t)\left(
\int_{\Gamma }^{t}1+\left\Vert u_{3}(\tau )\right\Vert _{L^{\alpha }}^{\frac{%
30\alpha }{7\alpha -45}}d\tau \right) ^{\frac{3\alpha -10}{4(\alpha -2)}%
}+C\left\Vert b_{3}(\Gamma )\right\Vert _{L^{10}}^{\frac{5}{2}}\epsilon ^{%
\frac{7}{8}}\mathcal{J}^{\frac{1}{4}}(t) \\
&&+C\epsilon ^{\frac{7}{8}}\mathcal{J}^{\frac{1}{4}}(t)\mathcal{L}^{\frac{5}{%
3}}(t)\left[ \left( \int_{\Gamma }^{t}\left\Vert \nabla b(\tau )\right\Vert
_{L^{2}}^{2}d\tau \right) ^{\frac{5\alpha }{3(4\alpha -15)}}\left(
\int_{\Gamma }^{t}\left\Vert u_{3}(\tau )\right\Vert _{L^{\alpha }}^{\frac{%
30\alpha }{7\alpha -45}}d\tau \right) ^{1-\frac{5\alpha }{3(4\alpha -15)}}%
\right] ^{\frac{4\alpha -15}{4\alpha }} \\
&\leq &C+C\epsilon ^{\frac{1}{4}}\mathcal{J}^{\frac{3}{2}}(t)+C\left\Vert
b_{3}(\Gamma )\right\Vert _{L^{10}}^{\frac{5}{2}}\epsilon ^{\frac{7}{8}}%
\mathcal{J}^{\frac{1}{4}}(t)+C\epsilon ^{\frac{7}{8}}\mathcal{J}^{\frac{1}{4}%
}(t)\mathcal{L}^{\frac{5}{3}}(t)\left[ \epsilon +\int_{\Gamma
}^{t}\left\Vert u_{3}(\tau )\right\Vert _{L^{\alpha }}^{\frac{30\alpha }{%
7\alpha -45}}d\tau \right] ^{\frac{4\alpha -15}{4\alpha }} \\
&\leq &C+C\epsilon ^{\frac{1}{4}}\mathcal{J}^{\frac{3}{2}}(t)+C\left\Vert
b_{3}(\Gamma )\right\Vert _{L^{10}}^{\frac{5}{2}}\epsilon ^{\frac{7}{8}}%
\mathcal{J}^{\frac{1}{4}}(t)+C\epsilon ^{\frac{7}{8}}\mathcal{J}^{\frac{1}{4}%
}(t)\mathcal{L}^{\frac{5}{3}}(t)\left[ 1+\int_{\Gamma }^{t}\left\Vert
u_{3}(\tau )\right\Vert _{L^{\alpha }}^{\frac{30\alpha }{7\alpha -45}}d\tau %
\right] ^{\frac{4\alpha -15}{4\alpha }} \\
&\leq &C+C\epsilon ^{\frac{1}{4}}\mathcal{J}^{\frac{3}{2}}(t)+C\left\Vert
b_{3}(\Gamma )\right\Vert _{L^{10}}^{\frac{5}{2}}\epsilon ^{\frac{7}{8}}%
\mathcal{J}^{\frac{1}{4}}(t)+C\epsilon ^{\frac{21}{4}}\mathcal{J}^{\frac{3}{2%
}}(t)+\frac{5}{6}\mathcal{L}^{2}(t)
\end{eqnarray*}%
which leads to%
\begin{equation}
\mathcal{L}^{2}(t)\leq C+C\epsilon ^{\frac{1}{4}}\mathcal{J}^{\frac{3}{2}%
}(t)+C\left\Vert b_{3}(\Gamma )\right\Vert _{L^{10}}^{\frac{5}{2}}\epsilon ^{%
\frac{7}{8}}\mathcal{J}^{\frac{1}{4}}(t)+C\epsilon ^{\frac{21}{4}}\mathcal{J}%
^{\frac{3}{2}}(t).  \label{eq25}
\end{equation}

Now, we will establish the bounds of $H^{1}-$norm of the velocity and
magnetic field. In order to do it, we multiply the first and second
equations of (\ref{eq1.1}) by $-\Delta u$ and $-\Delta b$, respectively, and
integrate them over $\mathbb{R}^{3}$ with respect to the spatial variable.
Then, integration by parts gives the following identity:
\begin{eqnarray}
&&\frac{1}{2}\frac{d}{dt}(\left\Vert \nabla u\right\Vert
_{L^{2}}^{2}+\left\Vert \nabla b\right\Vert _{L^{2}}^{2})+\left\Vert \Delta
u\right\Vert _{L^{2}}^{2}+\left\Vert \Delta b\right\Vert _{L^{2}}^{2}  \notag
\\
&=&-\sum\limits_{i,j,k=1}^{3}\int_{\mathbb{R}^{3}}\partial
_{k}u_{i}\partial _{i}u_{j}\partial
_{k}u_{j}dx+\sum\limits_{i,j,k=1}^{3}\int_{\mathbb{R}^{3}}\partial
_{k}b_{i}\partial _{i}b_{j}\partial
_{k}u_{j}dx-\sum\limits_{i,j,k=1}^{3}\int_{\mathbb{R}^{3}}\partial
_{k}u_{i}\partial _{i}b_{j}\partial _{k}b_{j}dx  \notag \\
&&+\sum\limits_{i,j,k=1}^{3}\int_{\mathbb{R}^{3}}\partial _{k}b_{i}\partial
_{i}u_{j}\partial _{k}b_{j}dx.  \label{eq12}
\end{eqnarray}%
Applying the divergence-free condition, $\nabla \cdot u=\nabla \cdot b=0$,
by using the H\"{o}lder inequality, the interpolation inequality, and (\ref%
{eq8}), it follows that
\begin{eqnarray*}
&&\frac{1}{2}\frac{d}{dt}(\left\Vert \nabla u\right\Vert
_{L^{2}}^{2}+\left\Vert \nabla b\right\Vert _{L^{2}}^{2})+\left\Vert \Delta
u\right\Vert _{L^{2}}^{2}+\left\Vert \Delta b\right\Vert _{L^{2}}^{2} \\
&\leq &C\int_{\mathbb{R}^{3}}(\left\vert \nabla _{h}u\right\vert +\left\vert
\nabla _{h}b\right\vert )(\left\vert \nabla u\right\vert ^{2}+\left\vert
\nabla b\right\vert ^{2})dx \\
&\leq &C(\left\Vert \nabla _{h}u\right\Vert _{L^{2}}+\left\Vert \nabla
_{h}b\right\Vert _{L^{2}})(\left\Vert \nabla u\right\Vert
_{L^{4}}^{2}+\left\Vert \nabla b\right\Vert _{L^{4}}^{2}) \\
&\leq &C(\left\Vert \nabla _{h}u\right\Vert _{L^{2}}+\left\Vert \nabla
_{h}b\right\Vert _{L^{2}})(\left\Vert \nabla u\right\Vert _{L^{2}}^{\frac{1}{%
2}}\left\Vert \nabla u\right\Vert _{L^{6}}^{\frac{3}{2}}+\left\Vert \nabla
b\right\Vert _{L^{2}}^{\frac{1}{2}}\left\Vert \nabla b\right\Vert _{L^{6}}^{%
\frac{3}{2}}) \\
&\leq &C(\left\Vert \nabla _{h}u\right\Vert _{L^{2}}+\left\Vert \nabla
_{h}b\right\Vert _{L^{2}})(\left\Vert \nabla u\right\Vert _{L^{2}}^{\frac{1}{%
2}}\left\Vert \nabla _{h}\nabla u\right\Vert _{L^{2}}\left\Vert \Delta
u\right\Vert _{L^{2}}^{\frac{1}{2}}+\left\Vert \nabla b\right\Vert _{L^{2}}^{%
\frac{1}{2}}\left\Vert \nabla _{h}\nabla b\right\Vert _{L^{2}}\left\Vert
\Delta b\right\Vert _{L^{2}}^{\frac{1}{2}}).
\end{eqnarray*}%
Integrating this last inequality in time, we deduce that for all $\tau \in
\lbrack \Gamma ,t]$%
\begin{eqnarray}
\mathcal{J}^{2}(t) &\leq &\left\Vert \nabla u(\Gamma )\right\Vert
_{L^{2}}^{2}+\left\Vert \nabla b(\Gamma )\right\Vert _{L^{2}}^{2}+C\underset{%
\tau \in \lbrack \Gamma ,t]}{\sup }(\left\Vert \nabla _{h}u(\tau
)\right\Vert _{L^{2}}+\left\Vert \nabla _{h}b(\tau )\right\Vert _{L^{2}})
\notag \\
&&\times \left( \int_{\Gamma }^{t}\left\Vert \nabla u(\tau )\right\Vert
_{L^{2}}^{2}d\tau \right) ^{\frac{1}{4}}\left( \int_{\Gamma }^{t}\left\Vert
\nabla \nabla _{h}u(\tau )\right\Vert _{L^{2}}^{2}d\tau \right) ^{\frac{1}{2}%
}\left( \int_{\Gamma }^{t}\left\Vert \Delta u(\tau )\right\Vert
_{L^{2}}^{2}d\tau \right) ^{\frac{1}{4}}  \notag \\
&&+C\underset{\tau \in \lbrack \Gamma ,t]}{\sup }(\left\Vert \nabla
_{h}u(\tau )\right\Vert _{L^{2}}+\left\Vert \nabla _{h}b(\tau )\right\Vert
_{L^{2}})  \notag \\
&&\times \left( \int_{\Gamma }^{t}\left\Vert \nabla b(\tau )\right\Vert
_{L^{2}}^{2}d\tau \right) ^{\frac{1}{4}}\left( \int_{\Gamma }^{t}\left\Vert
\nabla \nabla _{h}b(\tau )\right\Vert _{L^{2}}^{2}d\tau \right) ^{\frac{1}{2}%
}\left( \int_{\Gamma }^{t}\left\Vert \Delta b(\tau )\right\Vert
_{L^{2}}^{2}d\tau \right) ^{\frac{1}{4}}  \notag \\
&\leq &\left\Vert \nabla u(\Gamma )\right\Vert _{L^{2}}^{2}+\left\Vert
\nabla b(\Gamma )\right\Vert _{L^{2}}^{2}+2C\mathcal{L}(t)\epsilon ^{\frac{1%
}{4}}\mathcal{L}(t)\mathcal{J}^{\frac{1}{2}}(t)  \notag \\
&=&\left\Vert \nabla u(\Gamma )\right\Vert _{L^{2}}^{2}+\left\Vert \nabla
b(\Gamma )\right\Vert _{L^{2}}^{2}+C\epsilon ^{\frac{1}{4}}\mathcal{L}^{2}(t)%
\mathcal{J}^{\frac{1}{2}}(t).  \label{eq26}
\end{eqnarray}%
Inserting (\ref{eq25}) into (\ref{eq26}) and taking $\epsilon $ small
enough, then it is easy to see that for all $\Gamma \leq t<T^{\ast }$, there
holds
\begin{eqnarray*}
\mathcal{J}^{2}(t) &\leq &\left\Vert \nabla u(\Gamma )\right\Vert
_{L^{2}}^{2}+\left\Vert \nabla b(\Gamma )\right\Vert _{L^{2}}^{2}+C\epsilon
^{\frac{1}{4}}\mathcal{J}^{\frac{1}{2}}(t)+C\epsilon ^{\frac{1}{2}}\mathcal{J%
}^{2}(t) \\
&&+C\left\Vert b_{3}(\Gamma )\right\Vert _{L^{10}}^{\frac{5}{2}}\epsilon ^{%
\frac{9}{8}}\mathcal{J}^{\frac{3}{4}}(t)+C\epsilon ^{\frac{11}{2}}\mathcal{J}%
^{2}(t) \\
&<&\infty ,
\end{eqnarray*}%
which proves%
\begin{equation*}
\underset{\Gamma \leq t<T^{\ast }}{\sup }\left\Vert \nabla u(t)\right\Vert
_{L^{2}}^{2}+\left\Vert \nabla b(t)\right\Vert _{L^{2}}^{2}<\infty .
\end{equation*}%
This completes the proof of Theorem \ref{th1}.
\end{pf}

\section{Aknoledgements}

This work was done while the second author was visiting the Catania
University in Italy. He would like to thank the hospitality and support of
the University, where this work was completed. This research is partially
supported by P.R.I.N. 2019. The third author wish to thank the support of
"RUDN University Program 5-100".

\end{document}